\documentclass{article}
\usepackage{amsmath,amsthm,amsfonts,amscd,eucal}
\usepackage{epsfig}
\usepackage{subfigure}

\numberwithin{equation}{section}

\def\ca{{\mathcal A}}
\def\cb{{\mathcal B}}
\def\cc{{\mathcal C}}

\def\ch{{\mathcal H}}
\def\cai{{\mathcal I}}

\def\ck{{\mathcal K}}
\def\cl{{\mathcal L}}

\def\cp{{\mathcal P}}

\def\bn{{\mathbb N}}
\def\br{{\mathbb R}}

\def\a{\alpha}

\def\g{\gamma}        
\def\d{\delta}        
\def\eps{\varepsilon}

\def\th{\vartheta}

\def\l{\lambda}       \def\La{\Lambda}
\def\m{\mu}
\def\n{\nu}

\def\r{\rho}
\def\s{\sigma}       \def\S{\Sigma}
    
\def\t{\tau}

\def\f{\varphi}

\def\o{\omega}        

\def\ov{\overline}

\def\itm#1{\item{$(#1)$}}
\DeclareMathOperator{\Lim}{Lim}

\DeclareMathOperator{\supp}{supp}
\DeclareMathOperator{\Tr}{Tr}
\DeclareMathOperator{\diam}{diam}
\DeclareMathOperator{\vol}{vol}

\def\S{\Sigma}
\def\emp{\emptyset}
\def\geo{\text{geo}}
\def\subc{\underline{\d}}
\def\supc{\overline{\d}}
\def\subd{\underline{d}}
\def\supd{\overline{d}}
\def\inc{\uparrow}
\def\dec{\downarrow}
\def\B{\ov{B}}

\def\e#1{e^{#1}}

\def\ml{\m_{\lim}}

\newtheorem{Thm}{Theorem}[section]
\newtheorem{Cor}[Thm]{Corollary}
\newtheorem{Prop}[Thm]{Proposition}
\newtheorem{Lemma}[Thm]{Lemma}
\theoremstyle{definition}
\newtheorem{Dfn}[Thm]{Definition}
\newtheorem{exmp}[Thm]{Example}
\newtheorem{Assump}[Thm]{Assumption}
\theoremstyle{remark}
\newtheorem{rem}[Thm]{Remark} 
 
\begin{document}
\title{\huge  Dimensions and spectral triples \\
for fractals in $\br^{N}$}
\author{Daniele Guido, Tommaso Isola}
\date{}
\markboth{Spectral triples for fractals in $\br^{N}$}
{Spectral triples for fractals in $\br^{N}$}
\maketitle
\bigskip\bigskip\noindent
 Dipartimento di Matematica, Universit\`a di Roma ``Tor
Vergata'', I--00133 Roma, Italy.  E-mail: {\tt guido@mat.uniroma2.it, 
isola@mat.uniroma2.it}

 \begin{abstract}
    Two spectral triples are introduced for a class of fractals in
    $\br^{N}$.

    The definitions of noncommutative Hausdorff dimension and
    noncommutative tangential dimensions, as well as the corresponding
    Hausdorff and Hausdorff-Besicovitch functionals considered in
    \cite{GuIs9}, are studied for the mentioned fractals endowed with
    these spectral triples, showing in many cases their correspondence
    with classical objects.  In particular, for any limit fractal, the
    Hausdorff-Besicovitch functionals do not depend on the generalized
    limit $\o$.
 \end{abstract}

 \setcounter{section}{-1}

 \section{Introduction.}\label{sec:zeroth}
 
 In this paper we extend the analysis we made in \cite{GuIs9} to
 fractals in $\br^{N}$, more precisely we define spectral triples for
 a class of fractals and compare the  classical measures,
 dimensions and metrics with the measures, dimensions and metrics
 obtained from the spectral triple, in the framework of A. Connes' 
 noncommutative geometry \cite{Co}.
 
 The class of fractals we consider is the class of limit fractals, namely
 fractals which can be defined as Hausdorff limits of sequences of
 compact sets obtained via sequences of contraction maps.  This class
 contains the self-similar fractals and is contained in the wider class of
 random fractals \cite{random}.  On any limit fractal, the described
 limit procedure produces also a family of limit measures $\m_{\a}$,
 $\a>0$.  Among limit fractals, we consider in particular the
 translation fractals, namely those for which the generating
 similarities of a given level have the same similarity parameter.  It
 turns out that for translation fractals all the limit measures
 $\m_{\a}$ coincide.
 
 For limit fractals we introduce spectral triples which generalise the
 one considered by Connes in \cite{Co} for Cantor-like fractals,
 namely are based on an approximation of the fractal with sequences of
 pairs of points.  In the first spectral triple, the sequences consist
 of all descendants, via the generating similarities, of one (or
 finitely many) ancestral pair.  In the second triple, among the
 descendants of a single ancestor via the generating similarities, we
 consider all parent-child pairs.
 
 In both cases, when translation fractals are considered, we prove
 that the noncommutative Hausdorff dimension and tangential dimensions
 defined in \cite{GuIs9} coincide with their classical counterparts
 computed in \cite{GuIs11,GuIs13}.  Let us recall that the noncommutative
 tangential dimensions are the extreme points of the traceability
 interval, namely of the set of (singular) traceability exponents for
 the inverse modulus of the Dirac operator.  Therefore any of these
 exponents gives rise to a singular trace $\t_{\o}$ which in turn
 defines a trace on the algebra $\ca$ of the spectral triple, hence,
 by Riesz theorem, a measure on the fractal.  For translation fractals
 all these measures coincide with the limit measure.  In the case of
 the parent-child triple, an analogous result holds for any limit fractal,
 i.e. the measure coming from the traceability exponent $\a$ coincides
 with the limit measure $\m_{\a}$.  As a consequence, the measure
 generated by a singular trace $\t_{\o}$ is well defined, namely does
 not depend on the generalised limit procedure $\o$.
 
 Finally we study the distance on the fractal induced {\it {\`a la 
 Connes}} by the spectral triple. In the case of the parent-child 
 triple, the noncommutative distance is always equivalent to the 
 Euclidean distance, namely they induce the same topology. Then we 
 compare the noncommutative distance with the Euclidean geodesic 
 distance, namely with the distance defined in terms of rectifiable 
 curves contained in the fractal (when they exist). We prove that the 
 identity map from the fractal endowed with the geodesic distance to 
 the fractal endowed with the noncommutative distance is Lipschitz.  
 As a consequence, when the Euclidean distance and the geodesic 
 distance are bi-Lipschitz, this holds for the noncommutative 
 distance too.

 \section{Classical aspects}\label{SecClassic}

 We start this Section by recalling known results on self-similar 
 fractals, then we introduce the class of limit fractals and their 
 limit measures, and give some examples. We then introduce an open set 
 condition which allows us to characterise the limit measures on the 
 fractal (Theorem \ref{uniquemualpha}), and to compute them in case of 
 translation limit fractals, under a mild assumption (Theorem 
 \ref{thm:haus.dim}). Finally, we recall the notions of tangential 
 dimensions for metric spaces and measures from \cite{GuIs11,GuIs13}.
 
 \subsection{Preliminaries}\label{Prelim}
 The general reference for this Subsection is \cite{Falconer}.
 \par \noindent{\bf{Hausdorff measure and dimension.}} Let $(X,\r)$ be a
 metric space, and let $h:[0,\infty) \to [0,\infty)$ be non-decreasing
 and right-continuous, with $h(0)=0$.  When $E\subset X$, define, for
 any $\d>0$, $\ch^{h}_{\d}(E) := \inf \{ \sum_{i=1}^{\infty} h(\diam
 A_{i}) : \cup_{i} A_{i} \supset E, \diam A_{i} \leq \d \}$.  Then the
 {\it Hausdorff-Besicovitch (outer) measure} of $E$ is defined as
 $$                                                                               
   \ch^{h} (E) := \lim_{\d\to0}\ch^{h}_{\d}(E).                                     
 $$                                                                               
 If $h(t) = t^{\a}$, $\ch^{\a}$ is called {\it Hausdorff (outer) 
 measure} of order $\a>0$.

 The number
 $$                                                                               
 d_{H}(E) := \sup \{ \a>0 : \ch^{\a}(E) = +\infty \} = \inf \{ \a>0 : 
 \ch^{\a}(E) = 0 \}
 $$                                                                               
 is called {\it Hausdorff dimension} of $E$.

 \noindent{\bf{Selfsimilar fractals.}} Let $\{w_{j} \}_{j=1,\ldots,p}$
 be contracting similarities of $\br^{N}$, $i.e.$ there are
 $\l_{j}\in(0,1)$ such that $\| w_{j}(x) - w_{j}(y) \| = \l_{j}
 \|x-y\|$, $x,y\in \br^{N}$.  Denote by $\ck(\br^{N})$ the family of
 all non-empty compact subsets of $\br^{N}$, endowed with the
 Hausdorff metric, which turns it into a complete metric space.  Then
 $W: K\in\ck(\br^{N}) \to \cup_{j=1}^{p} w_{j}(K) \in \ck(\br^{N})$ is
 a contraction.
 \begin{Dfn}
	 The unique non-empty compact subset $F$ of $\br^{N}$ such that
	 $$
	 F = W(F) = \bigcup_{j=1}^{p} w_{j}(F)
	 $$ 
	 is called the {\it self-similar fractal} defined by $\{w_{j} 
	 \}_{j=1,\ldots,p}$.
 \end{Dfn}

 If we denote by $Prob_\ck(\br^{N})$ the set of probability measures 
 on $\br^{N}$ with compact support endowed with the Hutchinson metric, 
 $i.e.$ $d(\mu,\nu) := \sup \{ |\int fd\mu - \int fd\nu| : \|f\|_{Lip} 
 \leq 1 \}$, then the map
 $$
 \begin{matrix}
	 T : & Prob_\ck(\br^{N}) &\to &Prob_\ck(\br^{N})\cr 
	 &\mu&\mapsto&\sum_{j=1}^{p} \l_{j}^{s}\mu\circ w_{j}^{-1}
 \end{matrix} 
 $$
 is a contraction, where $s>0$ is the unique real number, called similarity
 dimension, satisfying $\sum_{j=1}^{p} \l_{j}^{s} =1$. We then observe that if
 $\m$ has support $K$, then $T\m$ has support $W(K)$. Since the sequence
 $W^n(K)$ is convergent, it turns out that it is bounded, namely there exists a
 compact set $K_0$ containing the supports of all the measures $T^n\m$. But on
 the space $Prob(K_0)$ the Hutchinson metric induces the weak$^*$ topology, and
 this space is compact in such topology, hence complete in the Hutchinson
 metric. Therefore there exists a fixed point of $T$ in $Prob_\ck(\br^{N})$,
 which is of course unique.

 \noindent{\bf{Open Set Condition.}} The similarities $\{w_{j}
 \}_{j=1,\ldots,p}$ are said to satisfy the open set condition if
 there is a non-empty bounded open set $V\subset \br^{N}$ such that
 $\cup_{j=1}^{p} w_{j}(V) \subset V$ and $w_{i}(V)\cap w_{j}(V)
 =\emptyset$, $i\neq j$.  In this case $d_{H}(F) = s$, the similarity
 dimension, and the Hausdorff measure $\ch^{s}$ is non-trivial on $F$. 
 Therefore $\ch^{s}|_{F}$ is the unique (up to a constant factor)
 Borel measure $\mu$, with compact support, such that $\mu(A) =
 \sum_{j=1}^{p} \l_{j}^{s} \mu(w_{j}^{-1}(A))$, for any Borel subset
 $A$ of $\br^{N}$.

 \subsection{Limit fractals.}\label{sec:Limit}
 
 Several generalisations of the class of self-similar fractals have 
 been studied.  Here we propose a new one, that we call the class of 
 limit fractals.  For its construction we need the following theorem.
 
 \begin{Thm}\label{gencontraction}
 	Let $(X,\r)$ be a complete metric space, $T_{n}:X\to X$ be such that 
 	there are $\l_{n}\in(0,1)$ for which $\r(T_{n}x,T_{n}y) \leq \l_{n} 
 	\r(x,y)$, for $x,y\in X$. Assume $\sum_{n=1}^{\infty} 
 	\prod_{j=1}^{n} \l_{j} <\infty$, and there is $x\in X$ such that 
 	$\sup_{n\in\bn} \r(T_{n}x,x) <\infty$. Then
	\itm{i} $\sup_{n\in\bn} \r(T_{n}y,y) <\infty$, for any $y\in X$, 
	\itm{ii}	$\lim_{n\to\infty} T_{1}\circ T_{2}\circ\cdots\circ 
	T_{n}x=x_{0}\in X$ for any $x\in X$.
 \end{Thm}
 \begin{proof}
 	$(i)$ $\r(T_{n}y,y) \leq \r(T_{n}y,T_{n}x) + \r(T_{n}x,x) + \r(x,y) 
 	\leq (1+\l_{n})\r(x,y) + \r(T_{n}x,x)$, so that $\sup_{n\in\bn} 
 	\r(T_{n}y,y) \leq 2\r(x,y) + \sup_{n\in\bn} \r(T_{n}x,x) <\infty$.\\
	$(ii)$ Set $M:= \sup_{n\in\bn} \r(T_{n}x,x) <\infty$, and $S_{n} 
	:= T_{1}\circ T_{2}\circ\cdots\circ T_{n}$, $n\in\bn$.  As 
	$\r(S_{n+1}x, S_{n}x) \leq \l_{1}\l_{2}\cdots\l_{n}\r(T_{n+1}x,x) 
	\leq M \l_{1}\l_{2}\cdots\l_{n}$, there follows, for any 
	$p\in\bn$, $\r(S_{n+p}x,S_{n}x)\leq \r(S_{n+p}x,S_{n+p-1}x) + 
	\ldots + \r(S_{n+1}x,S_{n}x) \leq M \sum_{k=n}^{n+p-1}\prod_{j=1}^{k} 
	\l_{k} \leq M \sum_{k=n}^{\infty}\prod_{j=1}^{k}\l_{k} \to 0$, as 
	$n\to\infty$, that is $\{ S_{n}x \}$ is Cauchy in $X$. Therefore 
	there is $x_{0}\in X$ such that $S_{n}x\to x_{0}$. \\
	Let us prove that $x_{0}$ is independent of $x$. Indeed, if $y\in X$, 
	then $\r(S_{n}x,S_{n}y) \leq \l_{1}\l_{2}\cdots \l_{n} \r(x,y) \to 
	0$, as $n\to\infty$, so that $S_{n}x$ and $S_{n}y$ have the same 
	limit.
 \end{proof}
 
 \begin{rem} 
	 A sufficient condition for $\sum_{n=1}^{\infty} \prod_{j=1}^{n} 
	 \l_{j} <\infty$ to hold is $$\sup_{n\in\bn} \l_{n}<1.$$
 \end{rem}
 
 We now describe the class of limit fractals.  Let $\{w_{nj} \}$,
 $n\in\bn$, $j=1,\ldots,p_{n}$, be contracting similarities of
 $\br^{N}$, with contraction parameter $\l_{nj}\in(0,1)$.  Set, for
 any $n\in\bn$, $\Sigma_{n} := \{\s:\{1,\ldots,n\}\to \bn : \s(k) \in
 \{1,\ldots,p_{k}\}, k=1,\ldots,n \}$, $\S:= \cup_{n\in\bn} \S_{n}$,
 $\S_{\infty} := \{\s:\bn\to \bn : \s(k) \in \{1,\ldots,p_{k}\},
 k\in\bn \}$, and write $w_{\s} := w_{1\s(1)}\circ w_{2\s(2)} \circ
 \cdots \circ w_{n\s(n)}$, $\l_{\s} := \l_{1\s(1)}\l_{2\s(2)} 
 \cdots \l_{n\s(n)}$, for any $\s\in\Sigma_{n}$.  Assume
 $\ov{\l}:=\sup_{n,j} \l_{nj} < 1$ and $\{ w_{\s}(x) : \s\in\Sigma \}$
 is bounded, for some (hence any) $x\in\br^{N}$.  Then, by Theorem
 \ref{gencontraction}, the sequence of maps $W_{n} : K\in\ck(\br^{N})
 \to \cup_{j=1}^{p_{n}} w_{nj}(K) \in \ck(\br^{N})$ is such that $\{
 W_{1}\circ W_{2}\circ\cdots\circ W_{n}(K) \}$ has a limit in
 $\ck(\br^{N})$, which is independent of $K\in\ck(\br^{N})$.
 
 \begin{Dfn} 
	 The unique compact set $F$ which is the limit of $\{ W_{1}\circ 
	 W_{2}\circ\cdots\circ W_{n}(K) \}_{n\in\bn}$ is called the limit 
	 fractal defined by $\{w_{nj} \}$.  In the particular case that 
	 $\l_{nj}=\l_{n}$, $j=1,\ldots,p_{n}$, $n\in\bn$, $F$ is called a 
	 translation (limit) fractal.  
 \end{Dfn}
 
 \begin{exmp} 
 As an example we mention some fractals considered in \cite{Hambly}. 
 They are constructed as follows.  At each step the sides of an
 equilateral triangle are divided in $q\in\bn$ equal parts, so as to
 obtain $q^{2}$ equal equilateral triangles, and then all downward
 pointing triangles are removed, so that $\frac{q(q+1)}{2}$ triangles
 are left.  
 The corresponding map $W$ can therefore be described as the map which 
 contracts the original triangle (or any of its subsets) by a factor 
 $1/q$, and then puts a copy of it in each of the upward pointing triangles. 
 Setting $q_{j}=2$ if $(k-1)(2k-1)<j\leq (2k-1)k$ and
 $q_{j}=3$ if $ k(2k-1)<j\leq k(2k+1)$, $k=1,2,\dots$, we get a
 translation fractal with dimensions given by (see Theorem \ref{thm:guis11})
 $$
 \subc=\frac{\log3}{\log2}<\subd=\supd=\frac{\log18}{\log6}
 <\supc=\frac{\log6}{\log3},
 $$
 where $\subc,\supc,\subd,\supd$ denote the lower tangential, the
 upper tangential, the lower local and the upper local dimensions. 
 The first four steps ($q=2,3,3,2$) of the procedure above are shown
 in Figure \ref{fig:Sierp}.

  \begin{figure}[ht]
     \centering
     \subfigure{
     \psfig{file=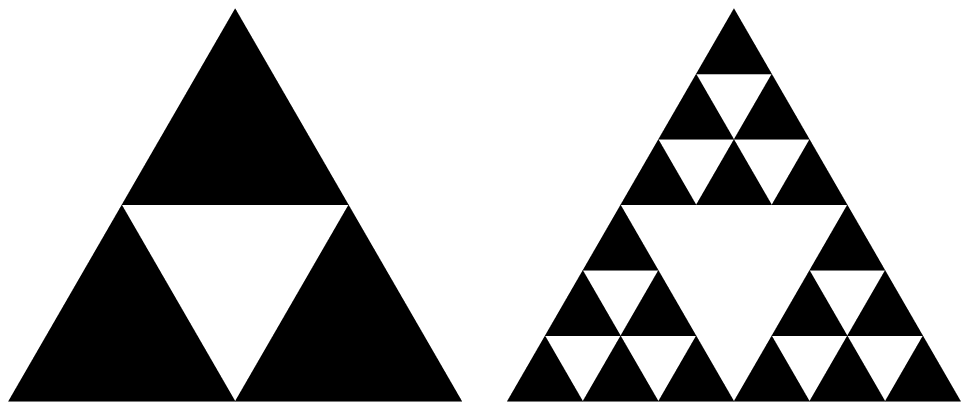,height=1.5in}}
     \hspace{0.3 in}
     \subfigure{
     \psfig{file=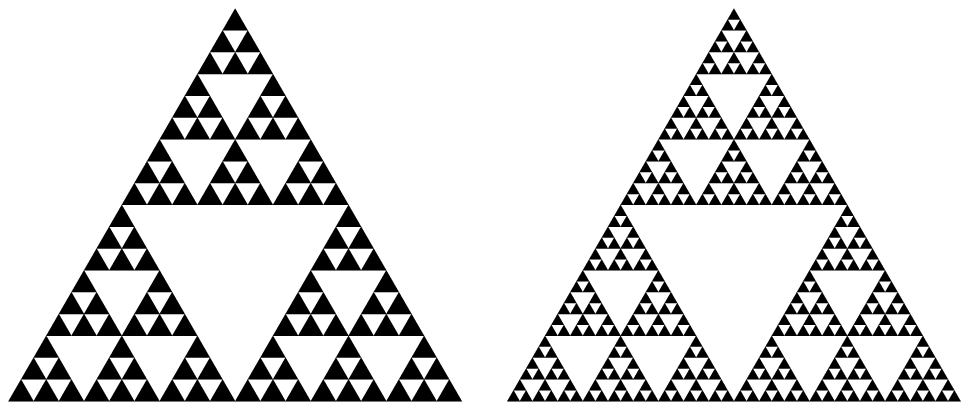,height=1.5in}}
     \caption{Modified Sierpinski}
     \label{fig:Sierp}
 \end{figure}

 The procedure considered above can, of course, be applied also to
 other shapes.  For example, at each step the sides of a square are
 divided in $2q+1$, $q\in\bn$, equal parts, so as to obtain
 $(2q+1)^{2}$ equal squares, and then $2q(q+1)$ squares are removed,
 so that to remain with a chessboard.  In particular, we may set
 $q_{j}=2$ if $k(2k+1)<j\leq (2k+1)(k+1)$ and $q_{j}=1$ if $
 k(2k-1)<j\leq k(2k+1)$, $k=0,1,2,\dots$, getting a translation
 fractal with dimensions given by (see Theorem \ref{thm:guis11})
 $$
 \subc=\frac{\log5}{\log3}<\subd=\supd=\frac{\log65}{\log15}
 <\supc=\frac{\log13}{\log5}.
 $$

 The first three steps ($q=1,2,1$) of this procedure are shown in Figure
 \ref{fig:Vics}.
 
 \begin{figure}[ht]
     \centering
     \subfigure{
     \psfig{file=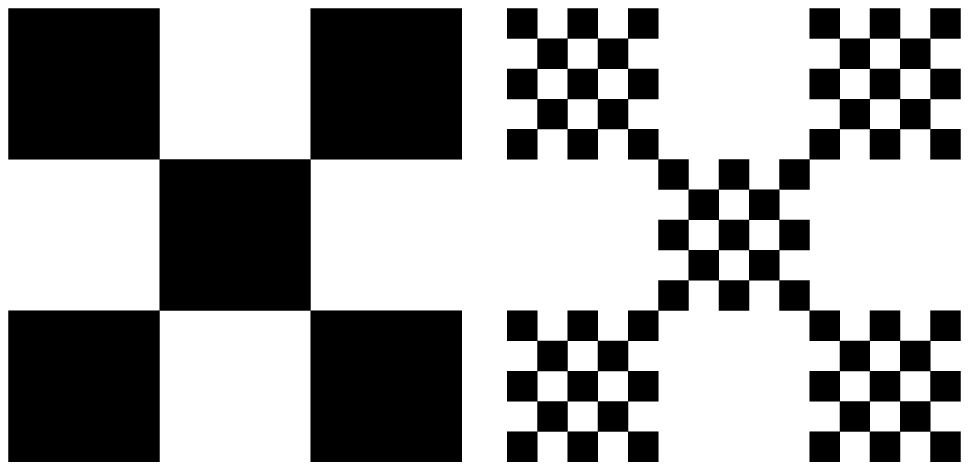,height=1.5in}}
     \hspace{0.3 in}
     \subfigure{
     \psfig{file=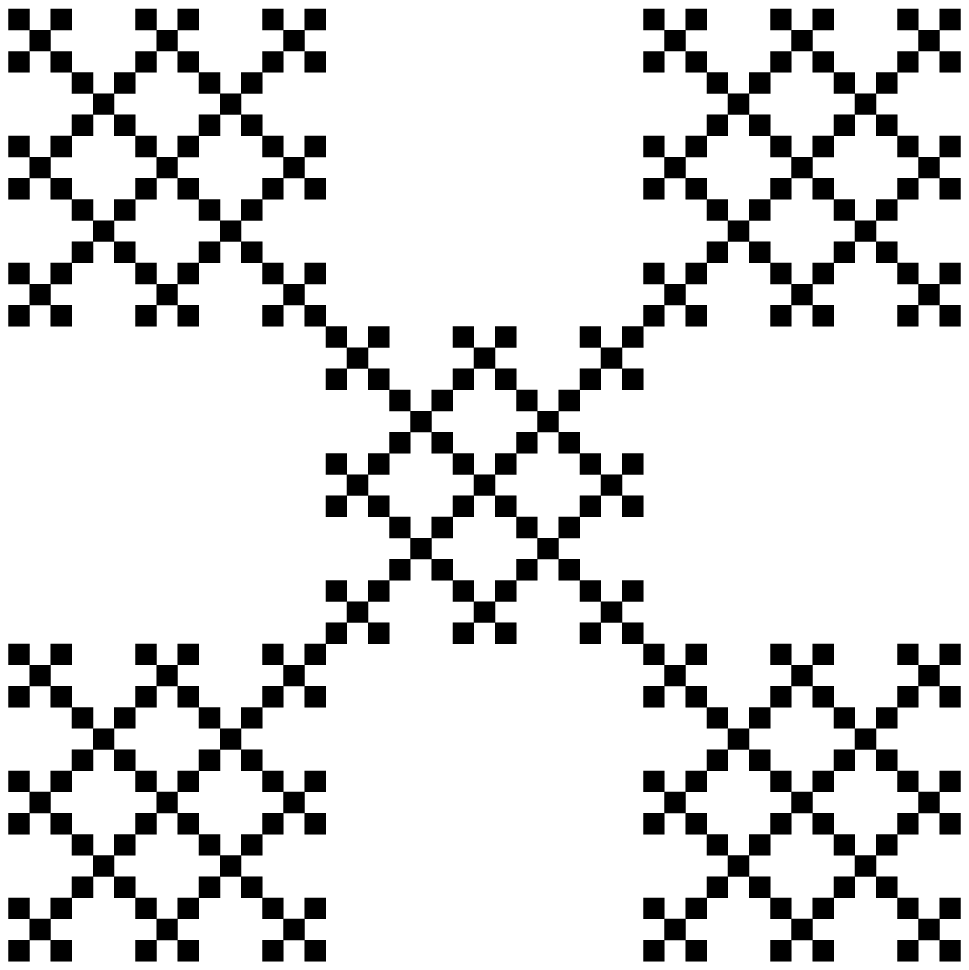,height=1.5in}}
     \caption{Modified Vicsek}
     \label{fig:Vics}
 \end{figure}

 \end{exmp}

 As before, we may consider the action of the similarities on 
 measures, besides that on sets.  Given $\a>0$ we set
 $$
 \begin{matrix}
	 T_n : & Prob_\ck(\br^{N}) &\to &Prob_\ck(\br^{N})\cr 
	 &\mu&\mapsto&\frac{1}{\sum_{j=1}^{p} \l_{nj}^{\a}} \sum_{j=1}^{p} 
	 \l_{nj}^{\a}\mu\circ w_{nj}^{-1}
 \end{matrix}
 $$
 and consider the sequence $\{T_{1}\circ T_{2}\circ\cdots\circ 
 T_{n}\m\}_{n\in\bn}$.  As before the supports of all such measures 
 are contained in a common compact set, therefore Theorem 
 \ref{gencontraction} applies and we get a unique limit measure 
 $\m_{\a}$, depending on the chosen $\a$.

As a consequence, if $\m_{0}$ is a probability measure and $\a>0$,
$\m_{n} := T_{1}\circ T_{2}\circ\cdots\circ T_{n}\m_{0}$ satisfies
\begin{equation}\label{1.1}
    \m_{n}(A) =\sum_{|\s|=n}c_{\s,\alpha}\m_{0}(w_{\s}^{-1}(A)),
\end{equation}
where	
\begin{equation}\label{1.2}
    c_{\s,\alpha} := \frac{\l_{\s}^{\a}}{\sum_{|\s'|=|\s|}\l_{\s'}^{\a}}.
\end{equation} 
In particular, if $F$ is a translation (limit) fractal, 
$$
c_{\s,\a}=\frac{1}{p_{1}p_{2}\cdots p_{|\s|}},\quad\forall\a>0,
$$
so that the limit measures $\m_{\a}$ all coincide, and will be denoted 
by $\ml$.

\subsection{Hausdorff dimension and limit measures}\label{sec:SetCond}

\begin{Assump} \label{ass:OSC}
    Open Set Condition: There is a bounded open set $V$ such that
    $w_{ni}(V)\subset V$ for any $n\in\bn,\, i\in\{1,\ldots,p_{n}\}$
    and $w_{n,i}V\cap w_{n,j}V=\emptyset$ if $i\ne j$.  We also assume
    that $V$ is regular, namely the Lebesgue measure of $V$ is equal
    to the Lebesgue measure of its closure $C$ and $V$ is equal to the
    interior of $C$.
\end{Assump}

OSC implies that $w_{ni} C\subset C$ for any $n\in\bn,\,
i\in\{1,\ldots,p_{n}\}$ and
$$
S_{N+1}C=W_1\cdots W_N\cdot W_{N+1}C\subset W_1\cdots W_N C=S_{N}C,
$$
namely $C_N=S_{N}C$ is a decreasing sequence of compact sets 
converging to $F$, $\cap C_N=F$.  Now set
\begin{align*}
V_{\s}&=w_{\s} V, \cr
C_{\s}&=w_{\s} C, \cr
F_{\s}&=w_{\s} C\cap F=C_{\s}\cap F. 
\end{align*}
Then 
\\
$(i)$
\begin{equation*}
\bigcup_{|\s|=N}F_{\s}=\bigcup_{|\s|=N}w_{\s} C\cap F=S_N C\cap F=C_N\cap 
F=F,
\end{equation*}
$(ii)$\quad if $\s\geq \s'$, namely $\s'$ is a truncation of $\s$, then
$V_{\s}\subset V_{\s'}$, \\
$(iii)$\quad if $\s,\s'$ are not ordered, $V_{\s}\cap V_{\s'}=\emptyset$.  

Moreover, in this case, $C_{\s}\cap V_{\s'}=\emptyset$.  In fact, if $x\in 
C_{\s}\cap V_{\s'}$, $x\in\partial V_{\s}$, hence there is a sequence $x_n\to 
x$, $x_n\in V_{\s}$, therefore $x_n$ is eventually in $V_{\s'}$, against the 
hypothesis.

\begin{Thm} \label{uniquemualpha}
    Let $F$ be a limit fractal satisfying OSC, with $\vol(V)=\vol(C)$. 
    Then the limit measure $\m_{\a}$ is the unique probability measure
    satisfying the following property: for any subset $\cai$ of 
    $\S_{n}$
    \begin{equation}\label{ineqCsa}
	\m_{\a}(V_{\cai}) \leq \sum_{\s\in\cai}c_{\s,\a} \leq 
	\m_{\a}(C_{\cai}),
    \end{equation}
    where we set $C_{\cai}=\cup_{\s\in\cai}C_{\s}$, $V_{\cai}$ equal 
    to the interior of $C_{\cai}$ relative to $C$, and $c_{\s,\a}$ is 
    defined in (\ref{1.2}).
\end{Thm}

\begin{proof}
    Let $\m_{n}$ converging to $\m_{\a}$ be as described in
    (\ref{1.1}).  Since $\m_{\a}$ does not depend on the starting
    measure, we may set $\m_{0}$ as the normalized Lebesgue measure on
    $V$.  Then $\m_{n}$ eventually satisfies
    $\m_{n}(V_{\s})=\m_{n}(C_{\s})=c_{\s,\a}$, therefore $\m_{n}$
    converges as a sequence of functionals on the vector space
    generated by continuous functions and step functions constant on
    the $V_{\s}$'s, giving rise to a positive bounded functional
    $\tilde\m$ on such space.  \\
    Then $\m_{\a}(V_{\cai})$ may be defined as the supremum of $\int
    f\ d \m_{\a}$ with support of $f$ contained in $V_{\cai}$, hence
    is majorised by $\tilde\m(V_{\cai})=\sum_{\s\in\cai} c_{\s,\a}$. 
    The second inequality of (\ref{ineqCsa}) is proved analogously. 
    Now we show that these inequalities determine $\m_{\a}$ uniquely. 
    Let $\m$ be a probability measure satisfying (\ref{ineqCsa}).  We
    observe that, for any continuous function $f$, $\int f\ d\m$ is
    well approximated by the lower Riemann sums with step functions
    constant on the $V_{\s}$'s, as soon as $|\s|$ is big enough, since
    $\diam(V_{\s}) \leq \diam(V) (\ov{\l})^{|\s|}$.  Let us now
    consider the set $\{ \min_{C_{\s}}f :\s\in\S_{n}\}$ and denote its
    elements by $f_{1},\dots f_{k}$ in increasing order.  Then
    define $\cai_{j}$ as the set of $\s\in\S_{n}$ such that $\min_{
    C_{\s}}f\geq f_{j}$, $C_{j}$ as the union of the $C_{\s}$ for
    $\s\in\cai_{j}$, $V_{j}$ as the interior of $C_{j}$.  Then
    \begin{align*}
	f_{1}+\sum_{i=2}^{k}(f_{i}-f_{i-1})\m(V_{i}) 
	& \leq f_{1} + \sum_{i=2}^{k} (f_{i}-f_{i-1}) 
	\left( \sum_{\s\in\cai_{i}} c_{\s,\a}\right) \\
	&\leq f_{1}+\sum_{i=2}^{k}(f_{i}-f_{i-1})\m(C_{i})
	\leq\int f\ d\m.
    \end{align*}	
    The second term may be rewritten as 
    $$
    \sum_{\s\in\S_{n}}\min_{ C_{\s}}f\ c_{\s\a},
    $$
    therefore 
    $$
    \int f\ d\m = \lim_{n\to\infty} \sum_{\s\in\S_{n}}
    \min_{ C_{\s}}f\ c_{\s\a},
    $$
    hence there is only one probability measure satisfying (\ref{ineqCsa}).
\end{proof}

Let now $F$ be a translation fractal ($\l_{n,i}$ independent of $i$),
and, to avoid triviality, assume $p_{n}\geq2$ for any $n\in\bn$.  The
OSC condition implies $\vol(S_{n-1}C) = \sum_{i=1}^{p_{n}}
\vol(w_{ni}S_{n-1}C) = p_{n}\l_{n}^{N}\vol(S_{n-1}C)$, so that
$2\l_{n}^{N}\leq1$, i.e. $\l_{n}\leq 2^{-1/N}$.  We set
\begin{equation}\label{LaP}
    \La_{n}=\prod_{i=1}^{n}\l_{i},\qquad P_{n}=\prod_{i=1}^{n}p_{i}.
\end{equation}
 
 \begin{Thm}\label{thm:haus.dim}
     Let $F$ be a translation fractal, with the notation above, and assume
     $p:=\sup_{n}p_{n}<\infty$.  Then 
    \begin{equation*}
	d_{H}(F) = \liminf_{n\to\infty} \frac{\log
	P_{n}}{\log 1/\La_{n}}.
    \end{equation*}
    Moreover the Hausdorff measure corresponding to $d:=d_{H}(F)$ is
    non trivial if and only if $\liminf(\log P_{n}-d \log 1/\La_{n})$
    is finite.
 \end{Thm} 
 \begin{proof}
    Let us consider the family $\cp$ of finite coverings of $F$, the
    subfamily $\cp(\S)$ of coverings made from sets of $\{C_{\s}:
    \s\in \S\}$, and the subfamily $\cp'(\S)$, whose coverings consist
    of $C_{\s}$, $\s\in \S$, $|\s|=const$.  If $P\in\cp$, $|P|$
    denotes the maximum diameter of the sets in $P$.  Clearly, for any
    $\a>0$, we have
    \begin{align*}
	\ch_{\a}(F) & = \lim_{\eps\to0}
	\inf_{\begin{smallmatrix}|P|\leq\eps\\
	P\in\cp\end{smallmatrix}} \sum_{E\in P}(\diam E)^{\a} \leq
	\lim_{\eps\to0} \inf_{\begin{smallmatrix}|P|\leq\eps\\
	P\in\cp(\S)\end{smallmatrix}} \sum_{E\in P}(\diam E)^{\a} \\
	& \leq \lim_{\eps\to0} \inf_{\begin{smallmatrix}|P|\leq\eps\\
        P\in\cp'(\S)\end{smallmatrix}} \sum_{E\in P}(\diam E)^{\a}.
    \end{align*}
    We shall show that the last two terms are indeed equal, and that
    the second term is majorised by a constant times the first, from
    which we derive 
    $$
    \ch_{\a}(F)\asymp\liminf_{n\to\infty}P_{n}\La_{n}^{\a},
    $$
    hence the required equality and the last statement.
    
    We may assume without restriction that the diameter of $V$
    is equal to one.  Then set $a:=\frac{\vol(V)}{\vol(B(0,2))}$,
    where $\vol$ denotes the Lebesgue measure.  Then the number of
    disjoint copies of $V$ intersecting a ball of radius $1$ is not
    greater than the number of disjoint copies of $V$ contained in a
    ball of radius $2$ which is in turn lower equal than $a^{-1}$.
	
    As a consequence, for any $x\in F$,
    \begin{equation}\label{VinBall}
	\#\{\s\in\S_{n}:\ C_{\s}\cap B_{F}(x,\La_{n})\ne\emptyset\}
	\leq a^{-1}.
    \end{equation}
    
    For any $\eps>0$, let $P=\bigcup_{i=1}^{n}E_{i}\in\cp$,
    $\diam E_{i}=r_{i}\leq\eps$.  Let now $x_{i}\in E_{i}$,
    $\La_{n_{i}+1}\leq r_{i}\leq\La_{n_{i}}$.
	
    Since the set $I_{n_{i}}\subset\S$ of multi-indices of length
    $n_{i}$ such that $\bigcup_{\s\in I_{n_{i}}}C_{\s}\supset
    B(x_{i},\La_{n_{i}})$ has cardinality majorised by $a^{-1}$, and
    any such $C_{\s}$ contains at most $p$ elements $C_{\s'}$, with
    $|\s'|=n_{i}+1$, then the set of multi-indices
    $I_{n_{i}+1}\subset\S$ of length $n_{i}+1$ such that
    $\bigcup_{\s\in I_{n_{i}+1}}C_{\s}\supset B_{F}(x_{i},\La_{n_{i}})$
    has cardinality majorised by $p/a$.  Therefore
    $$
    \bigcup_{i=1}^{n}\bigcup_{\s\in I_{n_{i}+1}}C_{\s}
    $$ 
    is a covering of $F$ of diameter less than $\eps$ and
    \begin{equation}
	\sum_{i=1}^{n}\sum_{\s\in I_{n_{i}+1}}(\diam C_{\s})^{\a} \leq
	\sum_{i=1}^{n}\frac{p}{a} \La_{n_{i}+1} ^{\a}
	\leq\frac{p}{a}\sum_{i=1}^{n}r_{i}^{\a}.
    \end{equation}
    As a consequence
    $$
    \lim_{\eps\to0} \inf_{\begin{smallmatrix}|P|\leq\eps\\
    P\in\cp(\S)\end{smallmatrix}} \sum_{E\in P}(\diam E)^{\a}
    \leq\frac{p}{a} \lim_{\eps\to0}\inf_{|P|\leq\eps}\sum_{E\in
    P}(\diam E)^{\a} = \frac{p}{a} \ch^{\a}(F).
    $$
    Now, for any $n_{1}\leq n_{0}$, let $P$ be the optimal covering of
    $F$ made of $C_{\s}$'s, with $n_{1}\leq|\s|\leq n_{0}$, namely
    minimizing $\sum_{\s\in P}(\diam C_{\s})^{\a}$, and choose $C_{\s_{0}}\in P$
    with $|\s_{0}|=n_{0}$.  This means that there is a $C_{\s}$,
    $|\s|=n_{0}-1$, which is optimally covered by some $C_{\s'}$'s of
    diameter $\La_{n_{0}}$.  Therefore this should be true for all
    other $\s$ of length $n_{0}$, namely the optimal covering is made
    of $C_{\s}$'s of the same size.  This shows the equality
    $$
    \lim_{\eps\to0} \inf_{\begin{smallmatrix}|P|\leq\eps\\
    P\in\cp(\S)\end{smallmatrix}} \sum_{E\in P}(\diam E)^{\a} =
    \lim_{\eps\to0} \inf_{\begin{smallmatrix}|P|\leq\eps\\
    P\in\cp'(\S)\end{smallmatrix}} \sum_{E\in P}(\diam E)^{\a}
    $$
    hence concludes the proof.
\end{proof}

\begin{rem}
	Let $F$ be a translation fractal, with the notation above, and
	assume $p:=\sup_{n}p_{n}<\infty$.  Let $G \subset F$ be closed. 
	Then, with $\cp(\S)$ denoting the family of finite coverings of
	$G$ made from sets in $\{V_{\s}:\s\in\S\}$, for any $\a>0$
	$$
	\ch^{\a}(G) \asymp \lim_{\eps\to0}
	\inf_{\begin{smallmatrix}|P|\leq\eps\\
	P\in\cp(\S)\end{smallmatrix}} \sum_{E\in P}(\diam E)^{\a}.
	$$
\end{rem}

It is clear that if the $F_{\s}$'s with $\s\in\S_{n}$ are essentially
disjoint w.r.t. $\m_{\a}$, then $\m_{\a}(F_{\s})=c_{\s,\a}$.  Now we
will discuss some conditions implying the vanishing of $\m_{\a}$ on
the intersections of the $F_{\s}$'s.  

\begin{Thm}\label{thm:special.transl}
	Let $F$ be a translation fractal for which $p:=\sup_{n}p_{n} <
	\infty$, $G$ a closed subset of $F$ s.t. $d_{H}(G) <
	d_{H}(F)$.  Then, $\ml(G)=0$.  \\
	As a consequence, if $d_{H}(w_{ni}C\cap w_{nj}C\cap
	F)<d_{H}(F)$, for any $n$, $i\ne j$, then
	$\ml(F_{\s})=\frac1{P_{|\s|}}$.
\end{Thm}
\begin{proof}
	Let $\a$ be s.t. $d_{H}(G) <\a < d_{H}(F)$, and $\eps>0$.  Then,
	from Theorem \ref{thm:haus.dim} and the Remark following it, there
	is $n_{0}\in\bn$,  s.t. $P_{n}\La_{n}^{\a}\geq 1$,
	for all $n\geq n_{0}$, and there is $\cai\subset\S$ s.t. $|\s|\geq
	n_{0}$, for all $\s\in\cai$, and $\sum_{\s\in\cai} \La_{\s}^{\a}
	\leq \eps$, and $G$ is contained in the interior of
	$\cup_{\s\in\cai} F_{\s}$.  By Urysohn's lemma, there is
	$f\in\cc(F)$, $0\leq f\leq1$, $f(x)=1,\, x\in G$, $\supp f\subset
	\cup_{\s\in\cai} F_{\s}$.  Then, with $\m_{k}$ as in (\ref{1.1}), 
	and $\m_{0}$ the normalised Lebesgue measure,
	\begin{align*}
		\ml(G) & \leq \int f\,d\m = \lim_{k\to\infty} \int
		f\,d\m_{k} \leq \lim_{k\to\infty}
		\m_{k}(\cup_{\s\in\cai} F_{\s}) \\
		&\leq \lim_{k\to\infty} \sum_{\s\in\cai} \m_{k}(F_{\s}) =
		\sum_{\s\in\cai} \frac{1}{P_{|\s|}} = \sum_{\s\in\cai}
		\frac{1}{P_{|\s|}\La_{|\s|}^{\a}} \La_{|\s|}^{\a} \\
		& \leq \sum_{\s\in\cai} \La_{|\s|}^{\a} \leq\eps. 
	\end{align*}
	The thesis follows.
\end{proof}

\subsection{Tangential dimensions}

 Let $(X,d)$ be a metric space, $E\subset X$.  Let us denote by
 $n(r,E) \equiv n_{r}(E)$, resp.  $\ov{n}(r,E) \equiv
 \overline{n}_{r}(E)$, the minimum number of open, resp.  closed,
 balls of radius $r$ necessary to cover $E$, and by $\n(r,E)
 \equiv \n_{r}(E)$ the maximum number of disjoint open balls of $E$ of radius
 $r$ contained in $E$.
 
 \begin{Dfn} {\rm \cite{GuIs11}}
     Let $(X,d)$ be a metric space, $E\subset X$, $x\in E$.  We call
     upper, resp.  lower tangential dimension of $E$ at $x$ the
     (possibly infinite) numbers
     \begin{align*}
	\underline{\d}_{E}(x) & := \liminf_{\l \to 0} \liminf_{r \to
	0} \frac{\log n(\l r, E\cap\B(x,r))}{\log 1/\l}, \\
	\overline{\d}_{E}(x) & := \limsup_{\l \to 0} \limsup_{r \to
	0} \frac{\log n(\l r, E\cap\B(x,r))}{\log 1/\l}.
    \end{align*}
 \end{Dfn}

 Tangential dimensions are invariant under bi-Lipschitz maps, and
 satisfy properties which are characteristic of a dimension function.

 Let $\m$ be a locally finite Borel measure, namely $\m$ is finite on
 bounded sets.
 
 Let us recall that the local dimensions of a measure at $x$ are
 defined as
 \begin{align*}
	 \underline{d}_{\m}(x)&=\liminf_{r\to0}
	 \frac{\log\m(B(x,r))}{\log r},\cr
	 \overline{d}_{\m}(x)&=\limsup_{r\to0}
	 \frac{\log\m(B(x,r))}{\log r}.
 \end{align*} 
  
 Now we introduce tangential dimensions for $\m$.
 
 \begin{Dfn}\label{meastgdims}{\rm \cite{GuIs13}}
     The lower and upper tangential dimensions of $\m$ are defined as
     \begin{align*}
	 \underline{\d}_{\m}(x) & := \liminf_{\l \to 0} \liminf_{r \to 0}
	 \frac{ \log \left( \frac{\m(B(x,r))} {\m(B(x,\l
	 r))} \right) } {\log 1/\l} \in [0,\infty], \\
	 \overline{\d}_{\m}(x) & := \limsup_{\l \to 0} \limsup_{r \to 0}
	 \frac{ \log \left( \frac{\m(B(x,r))}{\m(B(x,\l
	 r))} \right) } {\log 1/\l} \in [0,\infty].
     \end{align*}
 \end{Dfn}

  \begin{Thm}\label{eqdef}{\rm \cite{GuIs13}}
     Let $\m$ be a locally finite Borel measure on $X$. Then the following holds.
     $$
     \underline{\d}_{\m}(x) \leq \underline{d}_{\m}(x) \leq
     \overline{d}_{\m}(x) \leq \overline{\d}_{\m}(x).
     $$
\end{Thm}
 
Tangential dimensions are invariant under bi-Lipschitz maps.

 \begin{Thm}\label{thm:guis11}{\rm \cite{GuIs13}}
     Let $F$ be a translation fractal with the notations above,
     $\m=\ml$, and assume $p:=\sup_{n}p_{n}<\infty$.  Then
     \begin{align*}   
	\subd_{\m}(x) & = d_{H}(F) = \liminf_{n\to\infty} \frac{\log
	P_{n}}{\log 1/\La_{n}},\\
	\supd_{\m}(x) & = \limsup_{n\to\infty} \frac{\log P_{n}}{\log
	1/\La_{n}},\\
	\underline{\d}_{F}(x) &=\underline{\d}_{\m}(x) = \liminf_{n,k\to\infty}\frac{\log
	P_{n+k}-\log P_{n}} {\log 1/\La_{n+k}-\log 1/\La_{n}},\\
	\overline{\d}_{F}(x) &=\overline{\d}_{\m}(x) = \limsup_{n,k\to\infty}\frac{\log
	P_{n+k}-\log P_{n}} {\log 1/\La_{n+k}-\log 1/\La_{n}}.
    \end{align*}
 \end{Thm}

\section{Noncommutative aspects}\label{sec:Noncomm}

 \subsection{Singular traces on the compact operators of a Hilbert space.}

 In this section we recall the theory of singular traces on 
 $\cb(\ch)$ as it was developed by Dixmier \cite{Dixmier}, who first 
 showed their existence, and then in \cite{Varga}, \cite{AGPS} and 
 \cite{GuIs5}.

 A singular trace on $\cb(\ch)$ is a tracial weight vanishing on the
 finite rank projections.  Any tracial weight is finite on an ideal
 contained in $\ck(\ch)$ and may be decomposed as a sum of a singular
 trace and a multiple of the normal trace.  Therefore the study of
 (non-normal) traces on $\cb(\ch)$ is the same as the study of
 singular traces.  Moreover, making use of unitary invariance, a
 singular trace of a given operator should depend only on its
 eigenvalue asymptotics, namely, if $A$ and $B$ are positive compact
 operators on $\ch$ and $\m_n(A)=\m_n(B)+o(\m_n(B))$, $\m_n$ denoting
 the $n$-th eigenvalue, then $\t(A)=\t(B)$ for any singular trace
 $\t$.  The main problem about singular traces is therefore to detect
 which asymptotics may be ``resummed'' by a suitable singular trace,
 that is to say, which operators are singularly traceable.

 In order to state the most general result in this respect we need 
 some notation.  Let $A$ be a compact operator.  Then we denote by 
 $\{\m_n(A)\}$ the sequence of the eigenvalues of $|A|$, arranged in 
 non-increasing order and counted with multiplicity.  We consider also 
 the (integral) sequence $\{S_n(A)\}$ defined as follows:
 $$
 S_n(A)=
 \begin{cases}
	S^{\inc}_{n}(A) := \sum_{k=1}^{n} \m_k(A) & A\notin\cl^1 \cr 
	S^{\dec}_{n}(A) := \sum_{k=n+1}^{\infty}\m_k(A) & A\in\cl^1,
 \end{cases} 
 $$
 where $\cl^1$ denotes the ideal of trace-class operators.  We call a compact 
 operator {\it singularly traceable} if there exists a 
 singular trace which is finite non-zero on $|A|$.  We observe that 
 the domain of such singular trace should necessarily contain the 
 ideal $\cai(A)$ generated by $A$.  A compact operator is called {\it 
 eccentric} if
 \begin{equation}\label{eq:1.1}
	\frac{S_{2n_k}(A)}{S_{n_k}(A)}\to1
 \end{equation}
 for a suitable subsequence $n_k$.  Then the following theorem holds.

 \begin{Thm}\label{eccsingtrac} {\rm \cite{AGPS}}
     A positive compact operator $A$ is singularly traceable $iff$ it
     is eccentric.  In this case there exists a sequence $n_k$ such
     that both condition~(\ref{eq:1.1}) is satisfied and, for any
     generalised limit $\Lim_{\omega}$ on $\ell^\infty$, the positive
     functional
     $$
     \tau_\omega (B) =
     \begin{cases} 
	 \Lim_{\omega}\left(\left\{\frac{S_{n_k}(B)}{S_{n_k}(A)}\right\}\right)
	 &\quad B \in \cai(A)_+ \\
	 +\infty&\quad B \not\in \cai(A),\ B>0,
     \end{cases}
     $$
     is a singular trace whose domain is the ideal $\cai(A)$ generated
     by $A$.
 \end{Thm}

  The best known eigenvalue asymptotics giving rise to a singular trace 
 is $\m_{n}\sim\frac1n$, which implies $S_{n}\sim\log n$.  The corresponding 
 logarithmic singular trace is generally called Dixmier trace.

 \begin{Dfn} {\rm \cite{GuIs9}}
     If $A$ is a compact operator, set $f(t)=-\log\m_{A}(e^{t})$,
     $t\in\br$, where $\m_{A}$ is the extension of $\m_{n}(A)$ to a
     piecewise constant right continuous function on $[0,\infty)$. 
     Then define
     \begin{align*}
	 \underline{\d}(A) & := \left(\lim_{k\to\infty}
	 \limsup_{n\to\infty} \frac{\log
	 \frac{\m_{n}(A)}{\m_{kn}(A)}}{\log k}\right)^{-1} =
	 \left(\lim_{h\to\infty}\limsup_{t\to\infty}
	 \frac{f(t+h)-f(t)}{h}\right)^{-1}\\
	 \overline{\d}(A) & := \left(\lim_{k\to\infty}
	 \liminf_{n\to\infty} \frac{\log
	 \frac{\m_{n}(A)}{\m_{kn}(A)}}{\log k}\right)^{-1} =
	 \left(\lim_{h\to\infty}\liminf_{t\to\infty}
	 \frac{f(t+h)-f(t)}{h}\right)^{-1}\\
	 \overline{d}(A) & := \left(\liminf_{n\to\infty} \frac{\log
	 \m_{n}(A)}{\log 1/n}\right)^{-1}=
	 \left(\liminf_{t\to\infty}
	 \frac{f(t)}{t}\right)^{-1}.
     \end{align*} 
     Moreover, we say that $\a>0$ is an exponent of singular
     traceability for $A$ if $|A|^{\a}$ is singularly traceable.
 \end{Dfn}
 
 \begin{Thm}\label{Thm:SingTracExp} {\rm \cite{GuIs9}} 
     Let $A$ be a compact operator.  Then, the set of singular
     traceability exponents of $A$ is the relatively closed interval
     in $(0,\infty)$ whose endpoints are $\subc(A)$ and $\supc(A)$. 
     In particular, if $\supd(A)$ is finite nonzero, it is an exponent
     of singular traceability.
 \end{Thm}

 Note that the interval of singular traceability may be $(0,\infty)$,
 as shown in \cite{GuIs8}.  In \cite{GuIs12} the previous Theorem has
 been generalised to any semifinite factor, and some questions
 concerning the domain of a singular trace have been considered.

 \subsection{Singular traces and spectral triples}

 In this section we shall discuss some notions of dimension in 
 noncommutative geometry in the spirit of Hausdorff-Besicovitch 
 theory.

 As is known, the measure for a noncommutative manifold is defined 
 via a singular trace applied to a suitable power of some geometric 
 operator (e.g. the Dirac operator of the spectral triple of Alain 
 Connes).  Connes showed that such procedure recovers the usual volume 
 in the case of compact Riemannian manifolds, and more generally the 
 Hausdorff measure in some interesting examples \cite{Co}, Section 
 IV.3.

 Let us recall that $(\ca,\ch,D)$ is called a {\it spectral triple} when 
 $\ca$ is an algebra acting on the Hilbert space $\ch$, $D$ is a self 
 adjoint operator on the same Hilbert space such that $[D,a]$ is 
 bounded for any $a\in\ca$, and $D$ has compact resolvent.  In the 
 following we shall assume that $0$ is not an eigenvalue of $D$, the 
 general case being recovered by replacing $D$ with 
 $D|_{\ker(D)^\perp}$.  Such a triple is called $d^+$-summable, $d\in 
 (0,\infty)$, when $|D|^{-d}$ belongs to the Macaev ideal 
 $\cl^{1,\infty}=\{a:\frac{S_{n}^{\uparrow}(a)}{\log n}<\infty\}$.  
 
 The noncommutative version of the integral on functions is given by 
 the formula $\Tr_\omega(a|D|^{-d})$, where $\Tr_\omega$ is the 
 Dixmier trace, i.e. a singular trace summing logarithmic divergences.  
 By the arguments below, such integral can be non-trivial only if $d$ 
 is the Hausdorff dimension of the spectral triple, but even this 
 choice does not guarantee non-triviality.  However, if $d$ is finite 
 non-zero, we may always find a singular trace giving rise to a 
 non-trivial integral.

 \begin{Thm} \label{Thm:trace} {\rm \cite{GuIs9}} 
	 Let $(\ca,\ch,D)$ be a spectral triple.  If $s$ is an exponent of 
	 singular traceability for $|D|^{-1}$, namely there is a singular 
	 trace $\t$ which is non-trivial on the ideal generated by 
	 $|D|^{-s}$, then the functional $a\mapsto\t(a |D|^{-s})$ is a 
	 trace state (Hausdorff-Besicovitch functional) on the algebra 
	 $\ca$. 
 \end{Thm}
 
 \begin{rem}
	When $(\ca,\ch,D)$ is associated to an $n$-dimensional compact 
	manifold $M$, or to the fractal sets considered in \cite{Co}, the 
	singular trace is the Dixmier trace, and the associated functional 
	corresponds to the Hausdorff measure.  This fact, together with the 
	previous theorem, motivates the following definition.
 \end{rem}
 
 \begin{Dfn}\label{Dfn:dimensions} {\rm \cite{GuIs9}} 
	   Let $(\ca,\ch,D)$ be a spectral triple, $\Tr_{\o}$ the Dixmier 
	   trace.  
	   \itm{i} We call $\a$-dimensional Hausdorff functional 
	   the map $a\mapsto Tr_\omega(a |D|^{-\a})$; 
	   \itm{ii} we call (Hausdorff) dimension of the spectral
	   triple the number
	   $$
	   d(\ca,\ch,D) = \inf \{ d>0: |D|^{-d} \in \cl^{1,\infty}_{0} \} = 
	   \sup \{ d>0: |D|^{-d} \not\in \cl^{1,\infty} \},
	   $$
	   where $\cl^{1,\infty}_{0} = \{a : \frac{S_{n}^{\uparrow}
	   (a)} {\log n} \to 0 \}$.
	   \itm{iii} we call minimal, resp.  maximal dimension of the
	   spectral triple the quantity $\subc(|D|^{-1})$, resp. 
	   $\supc(|D|^{-1})$, hence $|D|^{\a}$ is singularly traceable
	   {\it iff} $\a\in[\subc,\supc]\cap(0,+\infty)$.
	   \itm{iv} For any $s$ between the minimal and the maximal
	   dimension, we call the corresponding trace state on the
	   algebra $\ca$ a Hausdorff-Besicovitch functional on
	   $(\ca,\ch,D)$.
 \end{Dfn}
  
 \begin{Thm}\label{unisingtrac} {\rm \cite{GuIs9}} 
     \itm{i} $d(\ca,\ch,D)=\supd(|D|^{-1})$. 
     \itm{ii} $d := d(\ca,\ch,D)$ is the unique exponent, if any, such
     that the $d$-dimensional Hausdorff functional is non-trivial.  
     \itm{iii} If $d\in (0,\infty)$, it is an exponent of singular
     traceability.
 \end{Thm}
  
 Let us observe that a singular trace, hence in particular the
 $\a$-dimensional Hausdorff functional, depends on a generalized limit
 procedure $\omega$, however its value is uniquely determined on the
 operators $a\in\ca$ such that $a|D|^{-d}$ is measurable in the sense
 of Connes \cite{Co}.  By an abuse of language we call measurable such
 operators.
 
 As in the commutative case, the dimension is the supremum of the 
 $\a$'s such that the $\a$-dimensional Hausdorff measure is everywhere 
 infinite and the infimum of the $\a$'s such that the $\a$-dimensional 
 Hausdorff measure is identically zero.  Concerning the non-triviality 
 of the $d$-dimensional Hausdorff functional, we have the same 
 situation as in the classical case.  Indeed, according to the 
 previous result, a non-trivial Hausdorff functional is unique (on 
 measurable operators) but does not necessarily exist.  In fact, if 
 the eigenvalue asymptotics of $D$ is e.g. $n\log n$, the Hausdorff 
 dimension is one, but the 1-dimensional Hausdorff measure gives the 
 null functional.
 
 However, if we consider all singular traces, not only the logarithmic 
 ones, and the corresponding trace functionals on $\ca$, as we said, 
 there exists a non trivial trace functional associated with 
 $d(\ca,\ch,D)\in(0,\infty)$, but $d(\ca,\ch,D)$ is not characterized 
 by this property.  In fact this is true if and only if the minimal 
 and the maximal dimension coincide.  A sufficient condition is the 
 following.
 
 \begin{Prop}\label{Prop:unique} {\rm \cite{GuIs9}} 
	 Let $(\ca,\ch,D)$ be a spectral triple with finite non-zero 
	 dimension $d$.  If there exists $\lim \frac{\m_n(D^{-1})} 
	 {\m_{2n}(D^{-1})} \in(1,\infty)$, $d$ is the unique exponent of 
	 singular traceability of $D^{-1}$.
 \end{Prop}

 \subsection{A spectral triple for fractals}\label{sec:Dirac}

 In this Subsection we introduce a spectral triple on limit 
 fractals, by extending an idea of Connes for Cantor-like fractals. We 
 compute its various dimensions, and recognise the noncommutative 
 Hausdorff functional as the one arising from the limit measure on the 
 fractal. Little can we say on the metric defined by this spectral 
 triple, so in the next Subsection we propose a different spectral 
 triple.
 
 Let $F$ be a limit fractal which satisfies OSC (se Assumption
 \ref{ass:OSC}) with respect to the open set $V$, and let $C$ be the 
 closure of $V$.  Choose two points $x,y\in C$, and denote with $r$
 their distance.  Also, construct the sequences $x_{\s}=w_{\s} x$,
 $y_{\s}=w_{\s} y$, ${\s}\in\S$ and note that
 $$
 \|x_{\s}-y_{\s}\|=r_{{\s}}:=r\l_{\s}.
 $$
 Then let $\ch$  be the $\ell^2$ space on the points $x_{\s}$, $y_{\s}$,
 and consider the natural representation of the Borel functions on $C$
 as multiplication operators on the elements of $\ch$.  Then let
 $$
 D := \bigoplus_{\s\in\S} 
 \frac{1}{\|x_{\s}-y_{\s}\|} 
     \begin{pmatrix}
		 0 & 1\\
		 1 & 0
	 \end{pmatrix}
 $$
 
 Now we consider the spectral triple $(\ca,\ch,D)$ where $\ch$ and $D$
 are defined as above, and $\ca$ is the algebra of continuous
 functions on $C$ such that $[D,f]$ is bounded.
 
 \begin{rem}\label{MoreAncestors} We may generalise the construction
 of the spectral triple by considering a finite number of ancestral pairs
 $\{x_{i},y_{i}\}$, e.g., for a Sierpinski like fractal as in Figure
 \ref{fig:Sierp}, the pairs of extreme points of the three sides of
 the original triangle.
 \end{rem}

 \begin{Thm}\label{thm:finadd.transl} 
     Let $F$ be a limit fractal, $(\ca,\ch,D)$ the spectral triple
     described above, $\a$ a singular traceability exponent for
     $|D|^{-1}$.  Then the Hausdorff-Besicovitch functional coincides,
     via Riesz Theorem, with the limit measure $\m_{\a}$.
 \end{Thm} 
 \begin{proof} 
     Let us first assume that the ancestral pair $\{x,y\}$ is 
     contained in $V$. Then the proof can be done as in Proposition 
     \ref{prop:finadd} and Theorem \ref{thm:finadd}.
     
     Now take a generic pair $\{x',y'\}$, with $\|x'-y'\|=r'$, and the
     corresponding spectral triple $(\ca,\ch,D')$, where the Hilbert
     spaces are naturally identified.  While the family of eigenvalues
     (with multiplicity) of $|D|^{-1}$ is given by
     $\{r\l_{\s}:\s\in\S\}$, each with multiplicity 2, the family of
     eigenvalues (with multiplicity) of $|D'|^{-1}$ is given by
     $\{r'\l_{\s}:\s\in\S\}$, each with multiplicity 2.  Therefore the
     spectral triples have the same set of traceability exponents, and
     if $\a$ is one of them, $\t$ is a singular trace such that
     $\t(|D|^{-\a})=1$, then $\t'(|D'|^{-\a})=1$, with
     $\t'=\left(\frac{r}{r'}\right)^{\a}\t$.
     
     Now let $f$ be a Lipschitz function on $C$, with Lipschitz
     constant $L$.  $f|D|^{-\a}$ is a multiplication operator, with
     eigenvalues (with multiplicity) $$\{r^{\a}\l_{\s}^{\a}f(x_{\s}),
     r^{\a}\l_{\s}^{\a}f(y_{\s}) : \s\in\S\}.$$ Then 
     $$
     f|D'|^{-\a}=\left(\frac{r'}{r}\right)^{\a}f|D|^{-\a}+R,
     $$ 
     where $R$ is a selfadjoint operator with eigenvalues (with
     multiplicity) $$\{(r')^{\a} \l_{\s}^{\a} (f(x'_{\s}) - f(x_{\s})),
     (r')^{\a} \l_{\s}^{\a} (f(y'_{\s}) - f(y_{\s})) : \s\in\S\}.$$
     Since $|f(x'_{\s}) - f(x_{\s})|\leq L\|x'_{\s} - x_{\s}\| =
     L\l_{\s} \|x' - x\|$, the operator $|R|$ can be majorised by a
     positive operator $S$ with eigenvalues $\{L M (r')^{\a}
     \l_{\s}^{\a+1} : \s\in\S\}$, each with multiplicity 2, where
     $M=\max(\|x' - x\| ,\|y' - y\|)$.  Clearly the operator $S$ is
     infinitesimal w.r.t. $|D|^{-\a}$, hence
     $$
     \t'(f|D'|^{-\a}) = \left(\frac{r'}{r}\right)^{\a} \t'(f|D|^{-\a})
     + \t'(R) = \t(f|D|^{-\a}) = \int f\ d\m_{\a}.
     $$
     Since Lipschitz functions are dense, we get the thesis.
 \end{proof}

 \begin{rem} 
     $(i)$ When $F$ is a translation fractal, and the Hausdorff
     dimension of the sets $F_{\s}\cap F_{\s'}$ is strictly lower than
     the Hausdorff dimension of $F$, $|\s|=|\s'|$, then Theorem
     \ref{thm:special.transl} applies, giving $$\ml(F_{\s}) =
     P_{|{\s}|}^{-1}.$$
     $(i)$ Let us observe that for self-similar fractals, there is
     only one exponent of singular traceability, namely the similarity
     dimension $s$, and the corresponding limit measure coincides with
     $\ch_{s}$.
 \end{rem}

 \begin{Thm}\label{Thm:spectrdim}
	Let $(\ca,D,\ch)$ be the spectral triple associated with a
	translation fractal $F$ where the similarities $w_{n,i}$,
	$i=1,\dots,p_{n}$ have scaling parameter $\l_{n}$.  Then the
	Hausdorff dimension is given by the formula
	$$
	d=\limsup_{n}\frac{\sum_1^n\log p_k}{\sum_1^n\log 1/\l_k}.
	$$
 \end{Thm}
 \begin{proof}
     The eigenvalues of $|D|^{-1}$ are given by $r$ with multiplicity
     $2$, $r\l_{1}$ with multiplicity $2p_{1}$, $r\l_{1}\l_{2}$ with
     multiplicity $2p_{1}p_{2}$, and so on.  Therefore, with
     $p_{0}:=1$, $\l_{0}:=1$, and $\La_{n}$, $P_{n}$ as in (\ref{LaP}),
     \begin{align*}
	 Tr(|D|^{-\a}) & = 2r^{\a} \sum_{k=0}^{\infty} \prod_{i=0}^{k}
	 (p_{k}\l_{k})^{\a} \\
	 & = 2r^{\a} \sum_{n=0}^{\infty} \e{\log P_{n}-\a\log 1/\La_{n}}
	 \\
	 & = 2r^{\a} \sum_{n=0}^{\infty} \exp\left(\log
	 P_{n}\left(1-\a\frac{\log 1/\La_{n}}{\log P_{n}}\right)\right).
     \end{align*}
     Denote by $d := \left(\liminf_{n\to\infty}
     \frac{\log 1/\La_{n}}{\log P_{n}} \right)^{-1}$.  Then, if $\a>d$,
     we get
     $$
     \liminf_{n\to\infty} \a\frac{\log 1/\La_{n}}{\log P_{n}} =
     \frac{\a}{d} >1
     $$
     which implies that, for any sufficiently small $\eps>0$ there is
     $n_{\eps}\in\bn$ such that, for all $n>n_{\eps}$,
     $$
     1 - \a\frac{\log 1/\La_{n}}{\log P_{n}} < -\eps.
     $$
     Since $P_{n}^{-\eps}\leq 2^{-n\eps}$, the series converges. 
     Whereas, if $\a<d$, as there is a subsequence $\{n_{k}\}$ such
     that $\frac{\log 1/\La_{n_{k}}}{\log P_{n_{k}}}\to d^{-1}$, we
     get
     $$
     1 - \a\frac{\log 1/\La_{n_{k}}}{\log P_{n_{k}}} \to 1-\frac{\a}{d} >0
     $$
     and the series diverges. Therefore $d(\ca,D,\ch) = d$. 
 \end{proof}
 
  \begin{Thm}\label{Thm:delta}
	Let $(\ca,\ch,D)$ be the spectral triple associated with a
	translation fractal $F$, where the
	similarities $w_{n,i}$, $i=1,\dots,p_{n}$ have scaling
	parameter $\l_{n}$.  Then
	\begin{align*}
		\subc(\ca,\ch,D)
		&=\liminf_{n,k}\frac{\sum_{j=n}^{n+k}\log p_{j}}
		{\sum_{j=n}^{n+k}\log 1/\l_j},\cr
		\supc(\ca,\ch,D)
		&=\limsup_{n,k}\frac{\sum_{j=n}^{n+k}\log p_{j}}
		{\sum_{j=n}^{n+k}\log 1/\l_j}.
	\end{align*}
 \end{Thm}
 \begin{proof}
     The eigenvalues of $|D|^{-1}$ are the numbers $r\La_{k}$, each
     with multiplicity $2P_{k}$.  Therefore, the quantity
     $\frac{1}{h}(\log 1/\m(e^{t+h})-\log1/\m(e^{t}))$, may be
     rewritten as
     \begin{equation}\label{ratioF}
	 \frac{\log 1/\La_{k}-\log 1/\La_{m}}
	 {\log\left(\sum_{j=0}^{k} P_{j}-\th_{k} P_{k}\right)
	 -\log\left(\sum_{j=0}^{m} P_{j}-\th'_{m} P_{m}\right)}
     \end{equation}
     for suitable constants $\th_{k}$, $\th'_{k}$ in $[0,1)$.  
	
     Let us observe that, since $p_{i}\geq2$, 
     \begin{align*}
	 \log\left(\sum_{j=0}^{k} P_{j}-\th_{k} P_{k}\right) - \log
	 P_{k} &\leq\log\frac{\sum_{j=0}^{k} P_{j}}{P_{k}}\\
	 &=\log\left(\sum_{j=0}^{k}\prod_{i=j+1}^{k}\frac1{p_{i}}
	 \right) \leq\log2.
     \end{align*}
     Since the denominator goes to infinity, additive perturbations of
     the numerator and of the denominator by bounded sequences do not
     alter the $\limsup$, resp.  $\liminf$, and the ratio above may be
     replaced by
     \begin{equation}\label{ratio2F}
	 \frac{\log 1/\La_{k}-\log 1/\La_{m}} {\log P_{k}-\log P_{m}}.
     \end{equation}
     Finally, since the denominator $\log P_{k}-\log P_{m}$ goes to
     infinity if and only if $k-m\to\infty$, the thesis follows.
 \end{proof}

 \begin{rem}
     As in Connes' book, we may introduce on $F$ the metric defined by
     $D$, namely
     $$
     d(x,y) := \sup \{ |f(x)-f(y)| : f\in\cc(F), \|[D,f]\|\leq 1 \}.
     $$
     However it is not true, in general, that $\iota : (F,\|\cdot\|)
     \to (F,d)$ is a homeomorphism.  For example, let $F$ be the
     Sierpinski gasket, $x,y\in F$ being two vertices of the
     enveloping triangle.  Then the metric $d$ gives value $+\infty$
     to any pair of points in $F$ sitting on a line which is not
     parallel to the side $x,y$.  Nevertheless, if we consider three
     ancestral pairs for $F$ as in remark \ref{MoreAncestors}, we get
     the Euclidean geodesic distance in $F$ (cf.  Theorem
     \ref{EuclidGeo}).
 \end{rem}

\subsection{A different spectral triple}\label{sec:Dirac2}

 In this final Subsection we construct a spectral triple on a large 
 subclass of limit fractals. We recognise the noncommutative Hausdorff 
 functionals as the ones arising from the limit measures, and compare 
 the ``noncommutative metric'' with the Euclidean geodesic distance. 
 In case of translation fractals, we compute the various dimensions 
 associated to the spectral triple.

 Let $F$ be a limit fractal satisfying OSC (see Assumption
 \ref{ass:OSC}).  Let $x_{ni}$ be the fixed point of $w_{ni}$, and
 assume that the set $W$ of all $x_{ni}$'s is not dense in $V$.  Fix
 $x_{\emp}\in V\setminus \ov{W}$, then there is $c>0$ s.t.
 $\|x_{\emp}-x_{ni}\| \geq c$, for all $n,i$, and
 \begin{align}
	\diam(C) & \geq \|x_{\emp} - w_{ni}x_{\emp}\| \geq \Big|
	\|x_{\emp} - x_{ni}\| - \|x_{ni} - w_{ni}x_{\emp}\| \Big| \notag\\
	& = (1-\l_{ni}) \|x_{\emp} - x_{ni}\| \geq (1-\ov{\l})c\label{eq:stime}.
 \end{align}

 Set $x_{\s}:=w_{\s}x_{\emp}$, and define
 $$
 D := \bigoplus_{\s\in\S}\bigoplus_{i=1}^{p_{|\s|+1}} 
 \frac{1}{\|x_{\s}-x_{\s\cdot i}\|} 
     \begin{pmatrix}
		 0 & 1\\
		 1 & 0
	 \end{pmatrix}
 $$
 acting on 
 $$\ch=\bigoplus_{\s\in\S}\bigoplus_{i=1}^{p_{|\s|+1}}
 \ell^{2}\{x_{\s},x_{\s\cdot i}\}$$
 
 In the following we shall consider the spectral triple $(\ca,\ch,D)$ 
 with $\ch$ and $D$ defined as above, and $\ca$ consisting of 
 continuous functions $f$ on $C$, acting on $\ch$ as
 \begin{align*}
 (f\xi)_{\s,i}(x_{\s})&=f(x_{\s})\xi_{\s,i}(x_{\s})
 \\
 (f\xi)_{\s,i}(x_{\s\cdot i})&=f(x_{\s\cdot i})\xi_{\s,i}(x_{\s\cdot i}),
 \end{align*}
 and for which $[D,f]$ is bounded.

 \subsubsection{Measures and dimensions}

 For $\a>0$, a singular traceability exponent of $|D|^{-1}$, set
 $$
 \int f d\m^{0}_{\a}:=\t(f|D|^{-\a})
 $$
 for any Borel function $f$, where $\t$ is a singular trace such that 
 $\t(|D|^{-\a})=1$.

\begin{Prop}\label{prop:finadd}
	With the above notation, $\m^{0}_{\a}(V_{\s}) =
	\m^{0}_{\a}(C_{\s}) = c_{\s,\a}$, for any $\s\in\S$.
\end{Prop}
\begin{proof}
	Let $\s,\s'\in\S_{n}$, then 
	$$\|x_{\s\cdot\r}-x_{\s\cdot\r\cdot i}\| = 
	\|w_{\s\cdot\r}x_{\emp}-w_{\s\cdot\r\cdot i}x_{\emp}\| = 
	\l_{\s} \|w_{\s}^{-1}w_{\s\cdot\r}x_{\emp} -
	w_{\s}^{-1}w_{\s\cdot\r\cdot i}x_{\emp}\|.$$
	As $w_{\s}^{-1}w_{\s\cdot\r}$ is independent of $\s$, the 
	operators $\l_{\s}^{-\a} \chi_{V_{\s}}|D|^{-\a}$ and 
	$\l_{\s'}^{-\a} \chi_{V_{\s'}}|D|^{-\a}$ have the same eigenvalues 
	(with multiplicity) up to finitely many. Therefore,
	$$
	1 = \t(|D|^{-\a}) = \sum_{\s'\in\S_{n}} \t(
	\chi_{V_{\s'}}|D|^{-\a} ) = \sum_{\s'\in\S_{n}}
	\frac{\l_{\s'}^{\a}}{\l_{\s}^{\a}}\t( \chi_{V_{\s}}|D|^{-\a}
	),
	$$
	so that $\m^{0}_{\a}(V_{\s}) = \t( \chi_{V_{\s}}|D|^{-\a} ) =
	c_{\s,\a}$. Since $\sum_{\s\in\S_{n}}c_{\s,\a}=1$, we get
	$\m^{0}_{\a}(C_{\s}) = c_{\s,\a}$.
\end{proof}

The set function $\m^{0}_{\a}$ is not $\s$-additive, however its
restriction to continuous functions gives rise to the
Hausdorff-Besicovitch functional on $\cc(F)=\ov\ca$.  The following
holds:

 \begin{Thm}\label{thm:finadd} 
      Let $\a$ be a traceability exponent.  Then the measure
      associated with the Hausdorff-Besicovitch functional via the
      Riesz Theorem, coincides with the limit measure $\m_{\a}$.  In
      particular it does not depend either on $x_{\emptyset}$ or on
      the generalised limit.
 \end{Thm} 
 \begin{proof} 
     Let us show that the regularization $\ov{\m}_{\a}$ of
     $\m^{0}_{\a}$ satisfies the inequalities (\ref{ineqCsa}).  Indeed
     $$
     \ov{\m}_{\a}(V_{\cai})\leq\m^{0}_{\a}(V_{\cai})\leq
     \sum_{\s\in\cai} \m^{0}_{\a}(C_{\s}) =\sum_{\s\in\cai}c_{\s,\a}.
     $$
     Moreover,
     $$
     \ov{\m}_{\a}(C_{\cai})\geq\m^{0}_{\a}(C_{\cai})\geq\sum_{\s\in\cai}
     \m^{0}_{\a}(V_{\s}) =\sum_{\s\in\cai}c_{\s,\a}.
     $$
     Then, by the uniqueness proved in Theorem \ref{uniquemualpha},
     $\ov{\m}_{\a}$ coincides with $\m_{\a}$.
 \end{proof}

 \begin{Thm}\label{specdim}
     Let $(\ca,D,\ch)$ be the spectral triple associated with a
     translation fractal $F$, where the similarities $w_{n,i}$,
     $i=1,\dots,p_{n}$ have scaling parameter $\l_{n}$, and 
     $\sup_{n}p_{n}<+\infty$.  Then 
     \begin{align*}
	 d(\ca,\ch,D)
	 &=\limsup_{n}\frac{\sum_1^n\log p_k}{\sum_1^n\log 1/\l_k}\\
	 \subc(\ca,\ch,D)
	 &=\liminf_{n,k}\frac{\sum_{j=n}^{n+k}\log p_{j}}
	 {\sum_{j=n}^{n+k}\log 1/\l_j},\\
	 \supc(\ca,\ch,D)
	 &=\limsup_{n,k}\frac{\sum_{j=n}^{n+k}\log p_{j}}
	 {\sum_{j=n}^{n+k}\log 1/\l_j}.
     \end{align*}
\end{Thm}
 \begin{proof}
     The eigenvalues of $|D|^{-1}$ are the numbers $\ov{\La}_{k,i} :=
     \La_{k}\|x_{\emp}-w_{k+1,i}x_{\emp}\|$, each with
     multiplicity $2 P_{k}$.  It follows from (\ref{eq:stime}) that
     there are $0<c_{1}<c_{2}$ s.t. $c_{1} \La_{k}\leq \ov{\La}_{k,i}
     \leq c_{2} \La_{k}$, for all $k\in\bn$. Therefore
     \begin{equation*}
	 Tr(|D|^{-\a}) 
	  = \sum_{k\in\bn} P_{k}\sum_{i=1}^{p_{k+1}} \ov{\La}_{k,i}^{\a}
	  \asymp \sum_{k} P_{k}\La_{k}^{\a},
     \end{equation*}
     and the first equality follows as in Theorem \ref{Thm:spectrdim}. 
     As for the others, the same computation in the proof of Theorem
     \ref{Thm:delta} can be performed.
 \end{proof}
 
 From Theorems \ref{thm:guis11}, \ref{specdim} we get
 
 \begin{Cor}
     Let $F$ be a translation (limit) fractal, $\m$ its limit measure, 
     $(\ca,\ch,D)$ its associated spectral triple. Then, for any $x\in 
     F$,
     \begin{align*}
     d(\ca,\ch,D) & = \supd_{\m}(x),\\
     \subc(\ca,\ch,D) & = \subc_{\m}(x),\\
     \supc(\ca,\ch,D) & = \supc_{\m}(x).
     \end{align*}
 \end{Cor}

 \subsubsection{Metrics}
 
 Let us now introduce on $F$ the metric defined by $D$ as in Connes' book,
 namely 
 $$
 d(x,y) := \sup \{ |f(x)-f(y)| : f\in\cc(C), \|[D,f]\|\leq 1 \}.
 $$
 
 \begin{Lemma}
     $\|x-y\|\leq d(x,y)$, for all $x,y\in F$.
 \end{Lemma}
 
 \begin{proof}
     Let $f$ be a Lipschitz function on $C$ with Lipschitz constant 
     $\leq1$. Then $\|[D,f]\|\leq1$, therefore $\|x-y\|\leq d(x,y)$.     
 \end{proof}

 Let us observe first that this distance and the Euclidean distance 
 induce the same topology.

 \begin{Lemma}\label{Lem:1}
     $(i)$ Let $\s\in\S_{n},\, x\in F_{\s}$.  Then
     $d(x,x_{\s}) \leq \frac{\diam(C)}{1-\overline{\l}}\l_{\s}$. 
     \item{$(ii)$} Let $\s,\s'\in\S_{n}$ be s.t. $F_{\s}\cap
     F_{\s'}\neq\emptyset$.  Then $d(x_{\s},x_{\s'}) \leq
     \frac{\diam(C)}{1-\overline{\l}} (\l_{\s}+\l_{\s'})$.
 \end{Lemma}
 \begin{proof}
     $(i)$. There is $\r\in\S_{\infty}$ s.t. $x=x_{\r}$ and the $n$-th 
     truncation $\r^{n}$ of $\r$ is equal to $\s$, $i.e.$, 
     $\r^{n}(i)=\s(i)$, $i=1,\ldots,n$.  Then
     \begin{align*}
	 d(x,x_{\s}) &\leq \sum_{k=|\s|}^{\infty}
	 d(x_{\r^{k}}, x_{\r^{k+1}}) \\
	 & \leq \sum_{k=|\s|}^{\infty} \l_{\r^{k}} \diam{(C)} \leq
	 \frac{\diam(C)}{1-\overline{\l}}\l_{\s} .
     \end{align*}	
     \\
     (ii) Let $x\in F_{\s}\cap F_{\s'}$.  Then $d(x_{\s},x_{\s'})
     \leq d(x,x_{\s}) + d(x,x_{\s'}) \leq
     \frac{\diam(C)}{1-\overline{\l}} (\l_{\s}+\l_{\s'})$.
 \end{proof}

 \begin{Thm}
     Let $F$ be a limit fractal. Then $\iota : (F,\|\cdot\|)\to (F,d)$ is 
     a homeomorphism.
 \end{Thm}
 \begin{proof}
     Assume $\iota$ is discontinuous in $x\in F$, so that there are 
     $c>0$, $\{x_{n}\}\subset F$ s.t. $\|x_{n} - x\|\to 0$, and 
     $d(x_{n},x)\geq c$, $n\in\bn$. From Lemma \ref{Lem:1} we obtain 
     $k\in\bn$ s.t., for $\s\in\S_{k}$, we have $\diam(F_{\s}) < 
     \frac{c}{2}$. As $\cup_{\s\in\S_{k}} F_{\s} = F$, there is 
     $\s\in\S_{k}$ s.t. $A:= \{n\in\bn : x_{n}\in F_{\s}\}$ is infinite, 
     hence $x\in F_{\s}$. Therefore $c \leq d(x_{n},x) \leq \diam(F_{\s}) < 
     \frac{c}{2}$, which is absurd.
 \end{proof}

The map $\iota$ is not bi-Lipschitz in general. This is true however 
when the Euclidean distance and the Euclidean geodesic distance in $F$ 
are bi-Lipschitz.

\begin{Lemma} \label{thm:necklace}
    Let $A\subset\br^{N}$ be a bounded open set, $\r\geq1$.  Then
    there are $k\in\bn$, $c>0$, such that, given similarities
    $\f_{1},\dots,\f_{k}$ such that $\f_{i}(A)\cap\f_{j}(A)=\emptyset$,
    $i\ne j$, with similarity parameters $\l_{1}\leq\dots\leq\l_{k}$
    satisfying $\l_{k}\leq\r\l_{1}$, and a rectifiable curve $\g$ such
    that $\g\cap\f_{i}(\ov{A}) \ne \emptyset$, $\forall i=1,\dots,k$, we have
    $\ell(\g)>c\l_{1}$.
\end{Lemma}
\begin{proof}
    Possibly rescaling, it is not restrictive to assume that
    $\l_{1}=1$.  Let $k\in\bn$ be s.t. the infimum length of a
    rectifiable curve intersecting the closure of $k$ dilated copies
    of $A$ with disjoint interior is $0$.  Then, for any $n\in\bn$, we
    get sequences of dilations $\f_{1n},\ldots,\f_{kn}$ such that
    $\f_{in}(A)\cap\f_{jn}(A)=\emptyset$, and a rectifiable curve
    $\g_{n}$ with $\ell(\g_{n}) < \frac{1}{n}$ such that
    $\g\cap\f_{in}(\ov{A}) \ne \emptyset$, $i=1,\dots,k$, $n\in\bn$. 
    Clearly it is not restrictive to assume that all curves $\g_{n}$
    start from the same point $x_{0}$, which implies that all the
    curves and copies of $A$ are contained in some compact set.  By
    the assumptions above, all dilations lie in a compact set.  Hence,
    possibly passing to a subsequence, we may assume that, for any
    $i=1,\dots k$, $\varphi_{in}$ converges to a dilation
    $\varphi_{i}$ of $\br^{N}$ and $\g_{n}\to\g$ in the Hausdorff
    topology.  Let us remark that in this way
    $\ov{A_{in}}\to\ov{A_{i}}$ in the Hausdorff topology, where
    $A_{i}:=\varphi_{i}(A)$, $A_{in}:=\varphi_{in}(A)$,
    $i=1,\ldots,k$, $n\in\bn$.  As $\diam(\g)=0$, $\g$ consists of a
    single point, indeed the point $x_{0}$.  Moreover
    $\g\cap\ov{A_{i}}\neq\emptyset$, for all $i$, i.e.
    $x_{0}\in\cap_{i=1}^{k}\ov{A_{i}}$.  \\
    We claim that $A_{1},\ldots,A_{k}$ are disjoint.  By
    contradiction, assume that $A_{i}\cap A_{j}$ is not empty, namely
    that there exist points $x_{i},x_{j}\in A$ with
    $\varphi_{i}(x_{i})=\varphi_{j}(x_{j})$, and let $r>0$ be s.t.
    $B(x_{i},r),B(x_{j},r) \subset A$.  Since the $\f_{i}$'s are
    dilations, this implies that $B(\varphi_{i}(x_{i}),r) \subset
    A_{i}$, and the same for $j$.  We have that $\varphi_{in}(x_{i})$
    and $\varphi_{jn}(x_{j})$ converge to the same point, hence, for a
    sufficiently large $n$, $\|\varphi_{in}(x_{i}) -
    \varphi_{jn}(x_{j})\| < r$, so that $A_{in}\cap
    A_{jn}\neq\emptyset$, which is absurd.  \\
    Since $x_{0}\in\cap_{i=1}^{k}\ov{A_{i}}$, so that $A_{i} \subset
    B(x_{0},\r\diam(A))$, we obtain $$\o_{N}\r^{N}(\diam(A))^{N} =
    \vol(B(x_{0},\r\diam(A))) \geq \sum_{i=1}^{k} \vol(A_{i}) \geq
    k\vol(A),$$
    with $\o_{N}$ the volume of the unit ball in $\br^{N}$. 
    Therefore, if we take the constant $k$ greater than
    $\frac{\o_{N}\r^{N}(\diam(A))^{N}}{\vol(A)}$, the infimum length
    of a rectifiable curve intersecting $k$ dilated copies of $A$ with
    disjoint interior cannot be $0$.
\end{proof}

\begin{Thm}\label{EuclidGeo}
    Let $F$ be a limit fractal in $\br^{N}$, $d_{\geo}$ the Euclidean
    geodesic distance in $F$, namely the distance defined in terms of
    rectifiable curves contained in the fractal (when they exist). 
    Assume $\underline{\l} := \inf_{n,i} \l_{ni} >0$.  Then there is
    $c\geq 1$ s.t. $\|x-y\| \leq d(x,y) \leq c d_{\geo}(x,y)$, $x,y\in
    F$.
\end{Thm}
\begin{proof}
    The first inequality was proved above, let us prove the second.  For any
    $\eps>0$, consider the set $\S(\eps)$ consisting of the
    multi-indices $\s$ such that $\l_{\s}\leq\eps$ but
    $\l_{\s^{n-1}}>\eps$, where $n=|\s|$ and $\s^{k}$ denotes the
    $k$-th truncation of $\s$ (as in the proof of Lemma \ref{Lem:1}). 
    Then $\l_{\s} = \l_{n,\s(n)} \l_{\s^{n-1}} > \underline{\l}
    \eps$.  It is clear that $\cup_{\s\in\S(\eps)}F_{\s}=F$, and
    $V_{\s}\cap V_{\s'}=\emptyset$ if $\s\ne\s'$ are in $\S(\eps)$.  \\
    Now let $x,y\in F$, $\eps>0$.  Choose a rectifiable curve $\g$ in
    $F$ connecting $x$ and $y$, with $\ell(\g) < 2 d_{geo}(x,y)$, and
    let $\s_{1},\dots\s_{k}$ be the elements of $\S(\eps)$ such that
    $C_{\s_{i}}\cap\g\ne\emptyset$, ordered in such a way that $x\in
    C_{\s_{1}}$, $y\in C_{\s_{k}}$, and $C_{\s_{i}}\cap C_{\s_{i+1}}
    \ne \emptyset$, $i=1,\dots,k-1$.
    \\
    By Lemma \ref{Lem:1} we get 
    \begin{align*}
	d(x,y)&\leq d(x,x_{\s_{1}}) + \sum_{i=1}^{k-1}
	d(x_{\s_{i}},x_{\s_{i+1}}) +d(x_{\s_{k}},y) \\
	&\leq \frac{2\diam(C)}{1-\ov{\l}}\sum_{i=1}^{k}\l_{\s_{i}} 
	\leq \frac{2\diam(C)k\eps}{1-\ov{\l}}.
    \end{align*}
    Let us notice that $k\to\infty$, when $\eps\to0$.

    Now let $k_{V}\in\bn$, $c_{V}>0$ be the constants associated to
    $V$ and to $1/\underline{\l}$ in Lemma \ref{thm:necklace}.  Then 
    \begin{align*}
	\ell(\g)& \geq [k/k_{V}]c_{V}\underline{\l}\eps  \\
	& \geq \frac{c_{V}\underline{\l}(1-\ov{\l})}{2\diam(C)} 
	\frac{[k/k_{V}]}{k}d(x,y).
    \end{align*}
    Passing to the limit for $\eps\to0$, we get $d(x,y)\leq
    \frac{4k_{V}\diam(C)} {c_{V}\underline{\l}(1-\ov{\l})}
    d_{\geo}(x,y)$, i.e. the thesis.
\end{proof}

\begin{rem}
    For many fractals, as Sierpinski gasket and carpet, Vicsek, 
    Lindstrom snowflake etc., $d_{\geo}$ and the Euclidean distance 
    are biLipschitz, hence also $d$ and the Euclidean distance are.
\end{rem}

 \end{document}